\input amstex

\def\BZ{{\Bbb Z}}
\def\CA{{\Cal A}}

\def\CC{{\Cal C}}
\def\CD{{\Cal D}}
\def\CF{{\Cal F}}
\def\CH{{\Cal H}}
\def\CK{{\Cal K}}

\def\CM{{\Cal M}}

\def\CO{{\Cal O}}

\def\CT{{\Cal T}}
\def\dpar{\partial}
\def\fg{{\frak g}}

\def\fm{{\frak m}}

\def\hCC{\widehat{{\Cal C}}}
\def\hCO{\widehat{{\Cal O}}}

\def\hS{\widehat S}
\def\hT{\widehat T}


\def\btu{\bigtriangleup}
\def\hAss{\widehat{Ass}}
\def\hCoass{\widehat{Coass}}
\def\hCocom{\widehat{Cocom}}
\def\hLie{\widehat{Lie}}
\def\hotimes{\widehat{\otimes}}
\def\hra{\hookrightarrow}
\def\hVect{\widehat{Vect}}
\def\iso{\buildrel\sim\over\longrightarrow}
 
\def\lla{\longleftarrow}
\def\lra{\longrightarrow}

\parskip=6pt

\documentstyle{amsppt}
\document
\magnification=1200
\NoBlackBoxes

\centerline{\bf REMARKS ON FORMAL DEFORMATIONS}
\bigskip
\centerline{\bf AND BATALIN-VILKOVISKY ALGEBRAS} 
\bigskip
\centerline{Vadim Schechtman}

\bigskip

\bigskip
\centerline{\bf Introduction}
\bigskip

This note consists of two parts. Part I is an exposition of 
(a part of) the V.Drinfeld's letter, [D]. The responsibility for 
the style lies entirely on the author of the present note. 

The sheaf of algebras of polyvector fields $\Lambda^\bullet\CT_X$ on a  
Calabi-Yau manifold $X$, equipped with the Schouten bracket,  
admits a structure of a Batalin-Vilkovisky algebra. This fact was probably 
first noticed by Z.Ran, [R]. Part II is devoted to some 
generalizations of this remark. 

To get the above BV structure, one uses a trivialization of the canonical 
bundle $\CK_X$. We note that, in fact, it is sufficient to have an 
integrable connection on $\CK_X$. Moreover,   
one has a canonical bijection 
between the set of integrable connections on $\CK_X$, or, what is the 
same, the set of {\it right} $\CD_X$-module structures on $\CO_X$, and 
the set of BV structures on $\Lambda^\bullet\CT_X$, cf. {\bf Theorem 4.3}. 
We give some generalizations of this correspondence, cf. \S 4.   

An interesting example of these data is a pair 
(a non compact CY manifold $X$, a function $f\in\Gamma(X,\CO_X)$). 
If we fix a trivialization $\CK_X\cong\CO_X$, we get the connection on 
$\CK_X$ defined by the exact form $df$, see \S 2. 

This note is an extended version of the talks given at Max-Planck-Institut 
f\"ur 
\newline Mathematik in January, 1998. I am very grateful to MPI for 
the hospitality and the excellent working conditions. 
I am grateful to Vladimir Hinich for the enlightening correspondence.    
I am especially obliged to Yu.I. Manin, 
whose question triggered this work. The contents of Sections 1 and 2 
of Part II was 
influenced by the M. Kontsevich's talk at MPI in January, 1997.

\newpage

\centerline{\bf TABLE OF CONTENTS}

\bigskip

Intoduction 1

\bigskip

{\bf Part I. Maurer-Cartan and formal deformations}

\bigskip

\S 1. Standard functors and Maurer-Cartan equations 3

\S 2. Koszul complex and Maurer-Cartan scheme 8

\S 3. Deformations 11

\bigskip

{\bf Part II. Milnor ring and Batalin-Vilkovisky algebras} 

\bigskip

\S 1. Koszul complex and Milnor ring 14

\S 2. Batalin-Vilkovisky algebras on Calabi-Yau manifolds 17

\S 3. Recollections on $D$-modules 21

\S 4. Connections and Batalin-Vilkovisky structures 24

\bigskip 

References 29

\newpage

\centerline{\bf PART I. MAURER-CARTAN AND FORMAL DEFORMATIONS}

\bigskip

\bigskip
\centerline{\bf \S 1. Standard functors and Maurer-Cartan equations}
\bigskip

{\bf 1.1.} In this Part, we fix once and for all a ground field $k$ 
of characteristic $0$. 
A "vector space" will mean a vector space over $k$; $\otimes$ will mean 
$\otimes_k$.  

$Vect$ will denote the category of vector spaces; it is a symmetric monoidal category, with respect to the tensor product. $Ass$ will denote the category 
of associative algebras in $Vect$ (without unit). The forgetful 
functor $Ass\lra Vect$ admits a left adjoint functor of a free 
associative algebra, to be denoted $\CF_{Ass}$. For $V\in Vect$, 
$\CF_{Ass}(V)$ is equal to the tensor algebra
$$
T^{\geq 1}(V)=\oplus_{n\geq 1}\ V^{\otimes n}, 
\eqno{(1.1)}
$$
with the obvious multiplication. 

Let $Com$ denote the category of commutative algebras in $Vect$ (assotiative, 
without unit). The forgetful functor $Com\lra Vect$ admits a left adjoint, 
$\CF_{Com}$. We have $\CF_{Com}(V)=S^{\geq 1}(V)$ where 
$$
S^{\geq 1}(V)=\oplus_{n\geq 1}\ S^n(V)
\eqno{(1.2)}
$$
Here $S^n(V)$ is the $n$-th symmetric power of $V$, i.e. the quotient 
of $V^{\otimes n}$ over the obvious action of the symmetric group 
$\Sigma_n$. 

The canonical projection $V^{\otimes n}\lra S^n(V)$ admits a canonical 
splitting $i:\ S^n(V)\hra V^{\otimes n}$, given by 
$$
i(x_1\cdot\ldots\cdot x_n)=\frac{1}{n!}\sum_{\sigma\in\Sigma_n}\ 
x_{\sigma(1)}\otimes\ldots\otimes x_{\sigma(n)}
\eqno{(1.3)}
$$
We denote by $Lie$ the category of Lie algebras in $Vect$, and by 
$\CF_{Lie}$ the functor of a free Lie algebra, the left adjoint 
to the forgetful functor $Lie \lra Vect$. 

{\bf 1.2.} Let $C$ be a coassociative coalgebra in $Vect$, with 
comultiplication $\Delta:\ C\lra C\otimes C$. For $n\geq 2$, define 
a map 
$$
\Delta^{(n)}:\ C\lra C^{\otimes n}
\eqno{(1.4)}
$$
by induction. We set $\Delta^{(2)}=\Delta$, and 
$$
\Delta^{(n+1)}=\bigl(id_{A^{\otimes (n-1)}}\otimes\Delta\bigr)\circ 
\Delta^{(n)}
\eqno{(1.5)}
$$
Let us consider the following condition 

{\bf (F)} {\it For each $x\in C$, there exists $n$ such that $\Delta^{(n)}(x)=
0$.}

We denote by $Coass$ the category of coassociative coalgebras in $Vect$ 
(without counit) satisfying {\bf (F)}. The forgetful functor 
$Coass\lra Vect$ admits a {\it right} adjoint, the functor 
of a {\it cofree} coalgebra, to be denoted $\CF_{Coass}$. 
For $V\in Vect$, $\CF_{Coass}(V)=T^{\geq 1}(V)$ as a vector space. 
The comultiplication is defined by 
$$
\Delta (v_1\otimes\ldots\otimes v_n)=\sum_{i=1}^{n-1}\ 
(v_1\otimes\ldots\otimes v_i)\otimes (v_{i+1}\otimes\ldots\otimes v_n)
\eqno{(1.6)}
$$ 
We denote by $Cocom$ the category of cocommutative coalgebras in $Vect$ 
(coassociative, without counit). The forgetful functor 
$Cocom\lra Vect$ admits a right adjoint, the functor of a cofree cocommutative 
coalgebra, to be denoted by $\CF_{Cocom}$. For $V\in Vect$, 
$\CF_{Cocom}(V)=S^{\geq 1}(V)$ as a vector space. The comultiplication is 
defined as the composition
$$
S^{\geq 1}(V)\hra T^{\geq 1}(V)\lra T^{\geq 1}(V)\otimes T^{\geq 1}(V)
\lra S^{\geq 1}(V)\otimes S^{\geq 1}(V)
\eqno{(1.7)}
$$
Here the first arrow is the canonical injection (1.3), the second one is 
the comultiplication in $T^{\geq 1}(V)$, and the third one is the canonical 
projection. 

{\bf 1.3.} By a {\it graded vector space}, we mean a collection of vector 
spaces $V^\bullet=\{V^i\}_{i\in\BZ}$, indexed by integers. Graded vector spaces 
form the category $GVect$, with obvious morphisms. It is a monoidal tensor 
category. For $V^\bullet,\ W^\bullet\in GVect$, their tensor product 
$V^\bullet\otimes W^\bullet$ is defined by 
$$
(V^\bullet\otimes W^\bullet)^i=\oplus_{p+q=i}\ V^p\otimes W^q
\eqno{(1.8)}
$$
The associativity constraint is the obvious one. The commutativity 
isomorphisms $s_{V^\bullet,W^\bullet}:\ V^\bullet\otimes W^\bullet\iso 
W^\bullet\otimes V^\bullet$ are given by 
$$
s_{V^\bullet,W^\bullet}(a\otimes b)=(-1)^{pq}b\otimes a,
\eqno{(1.9)}
$$
for $a\in V^p,\ b\in W^q$. We often use the notation $|a|=p$. 
 
We have  the {\it shift functors} $[i]:\ GVect\lra GVect\ (i\in\BZ)$, given by 
$(V^\bullet[i])^j=V^{i+j}$. By a {\it morphism of degree $i$} between 
$V^\bullet$ and $W^\bullet$ we mean a morphism $V^\bullet\lra W^\bullet[i]$ 
in $GVect$. We define $\CH om(V^\bullet,W^\bullet)\in GVect$ by 
$$
\CH om(V^\bullet,W^\bullet)^i=Hom_{GVect}(V^\bullet,W^\bullet)
\eqno{(1.10)}
$$ 

We define the categories $GAss, GCoass, GCom, GCocom$ and $GLie$ 
as the categories of associative algebras, coassociative coalgebras, 
commutative algebras, cocommutative coalgebras (all without unit or counit)  
and Lie algebras in 
$GVect$. {\it The coalgebras are assumed to satisfy the condition {\bf (F)}.} 
We regard $Vect$, etc. as a full subcategory of $GVect$, etc., consisting 
of objects living in degree $0$.  

We have the functors of (co)free (co)algebras $\CF_{GAss}, \CF_{GCoass}, 
\CF_{GCom}, \CF_{GCocom}$ and $\CF_{GLie}$, defined as in the 
case of non-graded algebras. 

If $A^\bullet\in GAss$ and $B^\bullet\in GCoass$ then 
$\CH om(B^\bullet,A^\bullet)$ admits a canonnical structure of a 
coassociative coalgebra. For two maps $f\in \CH om(B^\bullet,A^\bullet)^i,\ 
g\in\CH om(B^\bullet,A^\bullet)^j$, their product $f\cdot g$ is defined 
as a composition
$$
f\cdot g:\ B^\bullet\lra B^\bullet\otimes B^\bullet
\buildrel{f\otimes g}\over\lra A^\bullet[i]\otimes A^\bullet[j]
\lra A^\bullet[i+j]
\eqno{(1.11)}
$$
Here the first arrow is the comultiplication in $B^\bullet$, and the third one 
is the multiplication in $A^\bullet$. The second one is defined by 
$$
(f\otimes g)(b\otimes c)=(-1)^{j|b|}f(b)\otimes g(c)
\eqno{(1.12)}
$$
Similarly, if $\fg^\bullet\in GLie$ and $B^\bullet\in GCocom$, 
$\CH om(B^\bullet,\fg^\bullet)$ admits a canonical structure of a graded 
Lie algebra. 

If $A^\bullet\in GCom$ then $A^\bullet\otimes\fg^\bullet$ is canonically 
a graded Lie algebra, with the bracket defined by 
$$
[a\otimes x,b\otimes y]=(-1)^{|x||a|}ab\otimes [x,y]
\eqno{(1.13)}
$$

A {\it differential graded (dg)} vector space, or a {\it complex}, 
is a graded vector space $V^\bullet$ equipped with a differential of degree 
$1$, $d:\ V^\bullet\lra V^\bullet[1],\ d^2=0$. We can repeat all the discussion 
of this Subsection, adding the word "differential" to "graded". The 
corresponding categories are denoted by $DGVect, DGAss$, etc. 

{\bf 1.4.} Let $A$ be an associative algebra. Consider the cofree graded  
coalgebra $\CF_{GCoass}(A[1])$. There exists a unique differential $d$  
on this graded space, making it a dg coassociative coalgebra, whose degree 
$-2$ component
$$
d^{-2}:\ \CF_{GCoass}(A[1])^{-2}=A\otimes A\lra \CF_{GCoass}(A[1])^{-1}=A
$$
coincides with the multiplication $\mu$ in $A$. Note that the condition $d^2=0$ 
is equivalent to the associativity of $\mu$. 

Let us denote this dg coalgebra by $\CC_{Coass}(A)$. This way we get a 
functor $\CC_{Coass}:\ Ass\lra DGCoass$ which extends naturally to a 
functor 
$$
\CC_{Coass}:\ DGAss\lra DGCoass
\eqno{(1.14)}
$$

Conversely, let $B$ be a coassociative coalgebra. Consider the free graded 
algebra $\CF_{GAss}(B[-1])$. There exists a unique differential $d$ on this 
space, making it a dg associative algebra, whose degree $1$ component
$$
d^1:\ \CF_{GAss}(B[-1])^1=B\lra \CF_{GAss}(B[-1])^2=B\otimes B
\eqno{(1.15)}
$$
coincides with the comultiplication $\nu$ in $B$. The condition $d^2=0$ 
is equivalent to the coassociativity of $\nu$. 

Let us denote this dg algebra by $\CC_{Ass}(B)$. This way we get a functor 
$\CC_{Ass}:\ Coass\lra DGAss$, which extends naturally to a functor 
$$
\CC_{Ass}:\ DGCoass\lra DGAss
\eqno{(1.16)}
$$

{\bf 1.5.} Let $A^{\bullet}$ be a dg associative algebra. The {\bf 
Maurer-Cartan equation} in $A^{\bullet}$ is the equation
$$
da+a\cdot a=0,
\eqno{(1.17)}
$$
$a\in A^1$. The set of all $a\in A^1$ satisfying (1.7) will be denoted 
$MC(A^{\bullet})$.

{\bf 1.7.} Let $\hVect$ denote the category whose objects are topological 
vector spaces $V$, complete in a linear topology. Thus, the topology 
admits a base of neighbourhoods of $0$ consisting of subspaces $V^{(i)}$, 
and $V$ is the inverse limit of the spaces with discrete topology,  
$$
V=\lim_{\lla}\ V/V^{(i)}
\eqno{(1.19)}
$$
Morphisms are continuous linear maps. This is a monoidal tensor category, 
with the tensor product $\hotimes$ given by
$$
V\hotimes W=\lim_{\lla}\ V/V^{(i)}\otimes W/W^{(j)}
\eqno{(1.20)}
$$
We regard the tensor category $Vect$ as embedded in $\hVect$, as a full 
subcategory of spaces with the discrete topology. 

Let $\hCoass, \hCocom$, etc. denote the category of coassociative coalgebras, 
cocommutative coalgebras, etc., in $\hVect$. {\it We do not impose 
on the coalgebras the finiteness condition {\bf (F)}.} 

The forgetful functor 
$\hCoass\lra\hVect$, resp., $\hCocom\lra\hVect$, admits a right adjoint 
$$
\CF_{\hCoass}(V)=\hT^{\geq 1}(V):=\prod_{n\geq 1}\ V^{\otimes n},
\eqno{(1.21)}
$$
resp., 
$$
\CF_{\hCocom}(V)=\hS^{\geq 1}(V):=\prod_{n\geq 1}\ S^n(V)
\eqno{(1.22)}
$$
We define the category $G\hCoass, DG\hCoass$, etc., in the obvious way. 

{\bf 1.9.} We define the functor 
$$
\hCC_{Coass}:\ DG\hAss\lra DG\hCoass
\eqno{(1.23)}
$$
in the same way as (1.14). 

Given $A^\bullet\in DG\hAss$, a solution $a$ of (1.17) gives rise 
to a $0$-cocycle 
of the dg coalgebra $\hCC_{Coass}(A^\bullet)$, given by 
$$
\sum_{n=1}^\infty\ a^{\otimes n}\in\hT^{\geq 1}(A^1)\subset \hCC_{Coass}
(A^\bullet)^0
\eqno{(1.24)}
$$     

{\bf 1.10. Theorem.} {\it Let $A^{\bullet}\in DGAss$, $B^{\bullet}\in 
DGCoass$. We have canonically
$$
Hom_{DGAss}(\CC_{Ass}(B^{\bullet}),A^{\bullet})=
MC(\CH om(B^{\bullet},A^{\bullet}))=
$$
$$
=Hom_{DGCoass}(B^{\bullet},\CC_{Coass}(A^{\bullet}))
\eqno{(1.25)}
$$
In particular, the functor $\CC_{Ass}$ is left adjoint to $\CC_{Coass}$.} 
$\btu$

{\bf 1.11.} Let $\fg$ be a Lie algebra. Consider the cofree coalgebra $\CF_{GCocom}(\fg[1])$ on the graded space $\fg[1]$. There exists 
a unique differential 
$$
d:\ \CF_{GCocom}(\fg[1])\lra \CF_{GCocom}(\fg[1])[1]
$$
making $\CF_{GCocom}\fg[1]$ a dg cocommutative coalgebra, 
and whose degree $-2$ component 
$$
d^{-2}:\ \CF_{Cocom}(\fg[1])^2=\Lambda^2(\fg)\lra \CF_{Cocom}(\fg[1])^1=\fg
$$
coincides with the bracket in $\fg$. Let us denote this dg coalgebra 
$\CC_{Cocom}(\fg)$. This way we get a functor $\CC_{Cocom}:\ Lie\lra DGCocom$, which 
extends naturally to a functor 
$$
\CC_{Cocom}:\ DGLie\lra DGCocom
\eqno{(1.26)}
$$
Conversely, let $B$ be a cocommutative coalgebra. Consider the 
free graded Lie algebra $\CF_{Lie}(B[-1])$. There exists a unique differential 
on this space, making it a dg Lie algebra, whose  degree $1$ component 
$$
d^1:\ \CF_{Lie}(B[-1])^1=B\lra \CF_{Lie}(B[-1])^2=S^2(B)
$$
coincides with the comultiplication in $B$. Let us denote this dg 
Lie algebra by $\CC_{Lie}(B)$. This way we get a functor 
$\CC_{Lie}:\ Cocom\lra DGLie$ which extends naturally to a functor
$$
\CC_{Lie}:\ DGCocom\lra DGLie
\eqno{(1.27)}
$$
In a similar manner, one defines the functor 
$$
\hCC_{Cocom}:\ DG\hLie\lra DG\hCocom
\eqno{(1.28)}
$$

{\bf 1.12.} Let $\fg^{\bullet}$ be a dg Lie algebra. 
The {\bf Maurer-Cartan equation} in $\fg^{\bullet}$ is the equation
$$
da+\frac{1}{2}[a,a]=0, 
\eqno{(1.30)}
$$
$a\in\fg^1$. 
The set of all $a\in\fg^1$ satisfying (1.3) will be denoted by 
$MC(\fg^{\bullet})$. 

A solution $a$ of (1.30) gives rise to a $0$-cocycle of the dg coalgebra 
$\hCC_{Cocom}(\fg^{\bullet})$, given by 
$$
\sum_{n=1}^{\infty}\ \frac{a^n}{n!}\in\hS^{\geq 1}(\fg^1)
\subset \hCC_{Cocom}(\fg^{\bullet})
\eqno{(1.31)}
$$

{\bf 1.13. Theorem.} {\it Let $\fg^{\bullet}\in DGLie$, 
$B^{\bullet}\in DGCocom$. We have canonically
$$
Hom_{DGLie}(\CC_{Lie}(B^{\bullet}),\fg^{\bullet})=
MC(\CH om(B^{\bullet},\fg^{\bullet}))=
$$
$$
=Hom_{DGCocom}(B^{\bullet},\CC_{Cocom}(\fg^{\bullet}))
\eqno{(1.32)}
$$
In particular, the functor $\CC_{Lie}$ is left adjoint to $\CC_{Cocom}$.}

Assume for simplicity that $B$ is concentrated in degree $0$. 
By definition, a map of graded Lie algebras 
$\alpha: \CC_{Lie}(B)\lra\fg^{\bullet}$ 
is the same as a map of vector spaces $\alpha^1: B\lra \fg^1$; the map 
$\alpha$ is compatible with the differentials iff the map $\alpha^1$ 
satisfies the Maurer-Cartan equation. This proves the first equality in this case. The rest is proved similarly. $\btu$ 

\bigskip

\newpage
\centerline{\bf \S 2. Koszul complex and Maurer-Cartan scheme}
\bigskip

{\it A. KOSZUL}

\bigskip

{\bf 2.0.} In this section we will use (dg) algebras and coalgebras with 
unit (resp., counit). We will use the notation $(\cdot)^+$ for the objects 
obtained by formally adjoining a unit or a counit. We denote by 
$Com^+,\ Cocom^+$ the categories of unital algebras, resp., coalgebras.  

For a (dg) vector space $V$, we denote
$$
S^\bullet(V)=\oplus_{n\geq 0}\ S^n(V)
\eqno{(2.0)}
$$
It is the underlying vector space for $\CF_{Com}(V)^+$ and 
$\CF_{Cocom}(V)^+$.    

{\bf 2.1.} Let $\fg^\bullet$ be a graded Lie algebra. Assume that 
$\fg^1$ is finite dimensional. Then the Maurer-Cartan equation (1.30) 
defines a closed subscheme of $\fg^1$, considered as an affine space, 
to be called the {\bf Maurer-Cartan scheme} of $\fg^\bullet$, and 
denoted also by $MC(\fg^\bullet)$. It represents the functor 
$Com^+\lra Sets$ which assigns to a $k$-algebra 
$B$ the set    
$MC(\fg^\bullet\otimes B)$, cf. (1.3).  

More explicitely, set 
$$
A=S^\bullet(\fg^{1*})
\eqno{(2.1)}
$$
$\fg^{1*}$ denotes the dual space. 
We have $\fg^1=Spec(A)$. 

Let us choose a base $\{e_i\}$ of $\fg^{2*}$. Define the elements 
$$
f_i=(e_i\circ d,e_i\circ [\ ,\ ])\in \fg^{1*}\oplus S^2(\fg^{1*})
\subset A
\eqno{(2.2)}
$$
Here $d:\ \fg^1\lra \fg^2$ is the component of the differential, 
and $[\ ,\ ]:\ S^2\fg^1\lra \fg^2$ is the component of the bracket. 

We have 
$$
MC(\fg^\bullet)=Spec\bigl(A/(f_i)\bigr)
\eqno{(2.3)}
$$
Here we have denoted by $(f_i)$ the ideal generated by all $f_i$. 

{\bf 2.2.} Let us assume that $\fg^i=0$ for $i\neq 1, 2$. Such 
a dg Lie algebra is the same as the following set of data: 

two vector spaces $V=\fg^1,\ W=\fg^2$, a linear map $d:\ V\lra W$, 
and a symmetric map $b:\ S^2(V)\lra W$. 

We have 
$$
\CC_{Cocom}(\fg^\bullet)^+= 
S^\bullet(V)\otimes\Lambda^\bullet(W)
\eqno{(2.4)}
$$
as a graded coalgebra. Here $\Lambda^\bullet(W)=
\oplus_{n\geq 0}\ \Lambda^n(W)$ is the exterior algebra.   
In (2.4), $V$ has grading $0$, and $W$ has grading $1$. 

The differential acts as 
$$
0\lra S^\bullet(V)\lra S^\bullet(V)\otimes W
\lra S^\bullet(V)\otimes\Lambda^2(W)\lra\ldots 
\eqno{(2.5)}
$$
This complex is dual to a Koszul complex. To make a precise statement,  
assume that $V, W$ are finite 
dimensional. Set $A=S^\bullet(V^*)$, as in (2.1).    

Then the (restricted) 
dual of (2.5) is 
$$
\ldots\lra\Lambda^2(W)\otimes A\lra W\otimes A\lra A\lra 0
\eqno{(2.6)}
$$
If we choose a basis $\{e_i\}$ of the linear space $W^*$, we get the elements 
$f_i\in \fm_A$, as in (2.2). The complex (2.5) is nothing but the Koszul 
complex $K(A;(f_i))$, cf. II, 1.1.  

{\bf 2.3. Corollary.} {\it In the assumptions of} 2.2, {\it one has 
canonical isomorphism of commutative algebras
$$
\hCO_{MC(\fg^\bullet);0}=\bigl(H^0(\CC_{Cocom}(\fg^\bullet)^+)\bigr)^*
\eqno{(2.7)}
$$}
Here in the left hand side stays the completion of the  
local ring of the origin of the scheme $MC(\fg^\bullet)$. 

{\bf 2.4.} More generally,  
$$
\CC_{Cocom}(\fg^\bullet)^+=S^\bullet(\fg^{odd})\otimes\Lambda^\bullet
(\fg^{ev})
\eqno{(2.8)}
$$
for an arbitrary dg Lie algebra $\fg^\bullet$. Here 
$\fg^{odd}=\oplus_{i\in\BZ}\ \fg^{2i+1},\ \fg^{ev}=
\oplus_{i\in\BZ}\ \fg^{2i}$.  

Assume that $\fg^i=0$ for $i\leq 0$. Then  
$$
\CC_{Cocom}(\fg^\bullet)^+=\oplus_{p=1}^\infty\ 
\oplus_{i_1,\ldots,i_p
\geq 0}\ \bigl(S^{i_1}(\fg^1)\otimes\Lambda^{i_2}(\fg^2)\otimes\ldots
\otimes F^{i_p}(\fg^p)\bigr)
\eqno{(2.9)}
$$
where $F=S$ if $p$ is odd, and $F=\Lambda$ if $p$ is even. In this space 
the component $\fg^p$ has degree $p-1$. With the differential, this 
complex looks as 
$$
0\lra S^\bullet(\fg^1)\lra S^\bullet(\fg^1)\otimes \fg^2\lra 
S^\bullet(\fg^1)\otimes\bigl(\Lambda^2(\fg^2)\oplus\fg^3)\lra\ldots
\eqno{(2.10)}
$$
In other words, the very beginning depends only on $\fg^{\leq 2}$. 
Therefore, we get 

{\bf 2.5. Theorem.} {\it Let $\fg^\bullet$ be a dg Lie algebra sitting 
in degrees $\geq 1$. Assume that $\fg^1$ is finite dimensional. We have 
a canonical isomorphism of commutative algebras
$$
\hCO_{MC(\fg^\bullet);0}=\bigl(H^0(\CC_{Cocom}(\fg^\bullet)^+)\bigr)^*
\eqno{(2.11)}
$$}

\bigskip

{\it B. CHEVALLEY-EILENBERG}

\bigskip

{\bf 2.6.} A dg Lie algebra $\fg^\bullet$ such that $\fg^i=0$ for $i\neq 0, 1$, 
is the same as the following set of data: 

a Lie algebra $U=\fg^0$, a $U$-module $V=\fg^1$, and a linear map  
$d:\ U\lra V$ such that $d([u,u'])=u\cdot du'-u'\cdot du$, for all 
$u,u'\in U$.   

Here we have denoted by $u\cdot v$ the result of the action of an element 
$u\in U$ on $v\in V$. 

We want to describe the dg coalgebra $\CC_{Cocom}(\fg^\bullet)^+$. 
We have 
$$
\CC_{Cocom}(\fg^\bullet)^+=\Lambda^\bullet(U)\otimes S^\bullet(V)
\eqno{(2.12)}
$$
as a graded coalgebra, where we assign the degree $-1$ to $U$, and 
the degree $0$ to $V$. 

Consider the cocommutative counital coalgebra $B=S^\bullet(V)$. 
The action of the Lie algebra $U$ on $V$ induces the $U$-action on 
all symmetric powers 
$$
U\otimes S^n(V)\lra S^n(V)
$$
On the other hand, the differential $d:\ U\lra V$ induces the maps 
$$
U\otimes S^n(V)\buildrel{d\otimes id}\over\lra V\otimes S^n(V)\lra S^{n+1}(V),
\eqno{(2.13)} 
$$
the last arrow being the multiplication. Adding up, we get the maps 
$$
U\otimes S^n(V)\lra S^n(V)\oplus S^{n+1}(V)
$$
Adding over $n$, we get a map 
$$
U\otimes B\lra B
\eqno{(2.14)}
$$
which defines the action of the Lie algebra $U$ by coderivations on 
the coalgebra $B$. 

Geometrically, this action can be interpreted as follows. Assume that 
the space $V$ is finite dimensional. Consider the commutative 
algebra $A=S^\bullet(V^*)$, and regard $V$ as an affine scheme 
$Spec(A)$. The Lie algebra $U$ acts on $V$ by affine vector fields; namely, 
to an element $u\in U$, there corresponds the vector field $\tau_u$, 
whose value at a point $v\in V$ is equal to $du+u\cdot v$. Here we have 
identified the tangent space $T_{V;v}$ with $V$. Therefore, we get 
the action of $U$ on the functions
$$
U\otimes A\lra A
\eqno{(2.15)}
$$ 
Now, the action (2.10) is nothing but the dual to (2.11). 

Consider the homological Chevalley-Eilenberg complex of the Lie algebra 
$U$ with coefficients in the $B$, where the action is (2.10), 
$$
C_{homol}(U;B):\ \ldots\lra\Lambda^2(U)\otimes B\lra U\otimes B\lra B\lra 0
\eqno{(2.16)}
$$
Let us regard it as living in degrees $\leq 0$. 

{\bf 2.7. Lemma.} {\it The complex} (2.16) {\it coincides with 
$\CC_{Cocom}(\fg^\bullet)^+$, if we use the identification} (2.12). 

{\bf 2.8. Corollary.} {\it Assume that $V$ is finite dimensional. Then we have 
a canonical isomorphism of pro-artinian algebras, 
$$
\hCO_{\CO_V;0}^U=\bigl(H^0(\CC_{Cocom}(\fg^\bullet)^+)\bigr)^*
\eqno{(2.17)}
$$}

Here the superscript $(\cdot)^U$ means the subspace of $U$-invariants. 

{\bf 2.9.} Let us introduce a shorthand notation 
$$
C(\fg^\bullet)=\CC_{Cocom}(\fg^\bullet)^+
\eqno{(2.18)}
$$
Assume that $\fg^i=0$ for $i<0$. The dg Lie algebra $\fg^{\geq 1}$ is 
a dg module over the Lie algebra $\fg^0$, therefore its symmetric powers 
are also dg $\fg^0$-modules, i.e. we have the action maps   
$$
\fg^0\otimes S^n(\fg^\bullet)\lra S^n(\fg^\bullet)
\eqno{(2.19)}
$$
On the other hand, we have the maps  
$$
\fg^0\otimes S^n(\fg^1)\lra S^{n+1}(\fg^1)
\eqno{(2.20)}
$$
defined as in (2.13). Adding the maps (2.19) and (2.20) and summing up 
over $n\geq 0$ and $p\geq 1$, we get a map 
$$
\fg^0\otimes S^\bullet(\fg^{\geq 1}[1])\lra S^\bullet(\fg^{\geq 1}[1])
\eqno{(2.21)}
$$
which defines a structure of a $\fg^0$-module on 
$S^\bullet(\fg^{\geq 1}[1])=C(\fg^{\geq 1})$. 

{\bf 2.10. Theorem.} {\it The cocommutative dg coalgebra 
$C(\fg^{\bullet})$ is canonically 
isomorphic to the Chevalley-Eilenberg complex 
$C_{homol}(\fg^0;C(\fg^{\geq 1}))$, where the action of $\fg^0$ on 
$C(\fg^{\geq 1})$ is defined in} (2.21). 

The coalgebra structure on the Chevalley-Eilenberg 
complex comes from the cocommutative coalgebra structure on 
$C(\fg^{\geq 1})$. 

{\bf 2.11.} For a more general discussion, I refer the reader to the 
very nice paper [H].  

\bigskip

\bigskip
\centerline{\bf \S 3. Deformations}
\bigskip

{\bf 3.1.} Let $C^{\bullet}$ be a graded coalgebra. Recall that a 
{\it coderivation} of degree $d$ of $C^\bullet$ is a linear mapping 
$D:\ C^\bullet\lra C^\bullet [d]$ satisfying the {\it co-Leibnitz rule}
$$
\Delta(D(x))=\sum\ (D(a)\otimes b+(-1)^{d\cdot |a|}a\otimes D(b))
\eqno{(3.1)}
$$
if $\Delta(x)=\sum\ a\otimes b$. Here $\Delta$ is the comultiplication in 
$B^\bullet$. The coderivations form a fraded Lie algebra, to be denoted 
$Coder(C^\bullet)$.  
The bracket is given by  
$$
[D_1,D_2]=D_1D_2-(-1)^{|D_1||D_2|}D_2D_1
\eqno{(3.2)}
$$ 

Let $V$ be a vector space. Consider the cofree graded coalgebra 
$\CF_{GCoass}(V[-1])$. Set 
$$
\fg^{\bullet}_V=Coder(\CF_{GCoass}(V[-1]))
\eqno{(3.3)}
$$  
As a graded vector space, 
$$
\fg^{\bullet}_V=\oplus_{i\geq 0}\ Hom (V^{\otimes i},V), 
\eqno{(3.4)}
$$
where $\fg^i=Hom(V^{\otimes (i+1)},V)$. 

{\bf 3.2. Lemma.} {\it Suppose we are given a map 
$$
f:\ V\otimes V\lra V
$$
It is associative if and only if, when considered as an element of $\fg^1_V$, 
it satisfies the equation
$$
\frac{1}{2}[f,f]=0
$$}

Given such $f$, if we extend it to the coderivation of degree $1$ 
of 
$$
\CF_{Coass}\ (V[-1])=\oplus_{i\geq 1}\ V^{\otimes i}
$$
its component 
$$
f_3:\ V^{\otimes 3}\lra V^{\otimes 2}
$$
is given by
$$
f_3(a\otimes b\otimes c)=f(a\otimes b)\otimes c-a\otimes f(b\otimes c), 
$$
by the co-Leibnitz rule. Therefore the component of 
$\frac{1}{2}[f,f]=f\cdot f$ acting from $V^{\otimes 3}$ to $V$, sends 
$a\otimes b\otimes c$ to $f(f(a\otimes b)\otimes c)-f(a\otimes f(b\otimes c))$. 
The Lemma follows. $\btu$. 

{\bf 3.3.} Let $(A, f: A\otimes A\lra A)$ be an associative algebra. 
Let us define the dg Lie algebra $\fg_f^\bullet$, which 
is equal to 
$$
Coder (\CF_{GCoass}(A[1]))=\oplus_{i\geq 0}\ Hom (A^{\otimes i},A)
\eqno{(3.5)}
$$
as a graded Lie algebra, and with the differential $d_f=ad(f)$. 

Here we consider $f$ as an element of $\fg_f^1$, and by the previous Lemma, 
$[f,f]=0$. 

Now, let $f'=f+h:\ A\otimes A\lra A$ be another candidate for the 
multiplication. We have
$$
\frac{1}{2}[f+h,f+h]=[f,h]+\frac{1}{2}[h,h]=d_fh+\frac{1}{2}[h,h]
$$
Let $M_A$ denote the set of associative multiplications on $A$. 

More generally, define $M_A$ as a functor $Com^+\lra Sets$ which assigns 
to a commutative unital $k$-algebra $B$ the set of associative 
multiplications on $A\otimes B$.  

We have proven 

{\bf 3.4. Lemma.} {\it Assume that $A$ is finite dimensional. Then 
$M_A$ is represented by an affine scheme.  
We have a canonical isomorphism of schemes 
$$
M_A=MC(\fg^\bullet_f)
$$
which takes $f\in M_A(k)$ to $0\in MC(\fg^\bullet_f)(k)$.} 

$\btu$ 

Applying Theorem 2.5, we get

{\bf 3.5. Theorem.} {\it We have a canonical isomorphism of formal 
commutative algebras
$$
\hCO_{M_A;f}=\bigl(H^0(\CC_{Cocom}(\fg^{\geq 1}_f)^+)\bigr)^*
\eqno{(3.6)}
$$}

\newpage

\centerline{\bf PART II. MILNOR RING AND BATALIN-VILKOVISKY ALGEBRAS}

\bigskip

\bigskip
\centerline{\bf \S 1. Koszul complex and Milnor ring}
\bigskip

{\bf 1.1. Koszul complex.} Let $A$ be a commutative ring, 
and $f_1,\ldots,f_r$ a sequence of its elements. 

We denote by 
$K^\bullet(A;f_1,\ldots,f_r)$ the complex
$$
0\lra\Lambda^r(A^r)\lra\Lambda^{r-1}(A^r)\lra\ldots\lra A^r\lra A\lra 0
\eqno{(1.1)}
$$
The exterior powers are over $A$. The differential is given by
$$
d(a_{i_1}\wedge\ldots\wedge a_{i_p})=\sum_{j=1}^p\ (-1)^{j-1}f_{i_j}a_{i_1}\wedge\ldots
\wedge\hat{a}_{i_j}\wedge\ldots\wedge a_{i_p}
\eqno{(1.2)}
$$
Here for $a\in A$ and $1\leq i\leq r$, we denote by $a_i\in A^r$ the element 
$(0,\ldots,a,\ldots,0)$, with $a$ on $i$-th place. 
We consider this complex as living in degrees from $-r$ to $0$.  

We have 
$$
K^\bullet(A;f_1,\ldots,f_r)=K^\bullet(A;f_1)\otimes_A\ldots\otimes_A 
K^\bullet(A;f_r)
\eqno{(1.3)}
$$
Equipped with the obvious wedge product, $K^\bullet(A;f_1,\ldots,f_r)$ 
becomes a commutative dg $A$-algebra. 

In fact, the differential (1.2) is the unique differential on the 
exterior algebra $\Lambda^\bullet(A^r)$,  
whose component 
$A^r\lra A$ is given by
$$
d(a_1,\ldots,a_r)=\sum_{i=1}^r\ f_ia_i,
\eqno{(1.4)}
$$
and making it a commutative dg algebra. 

{\bf 1.2.} Assume that $(f_1,\ldots,f_r)$ is a regular sequence, i.e. 
the operator of multiplication by $f_i$ on $A/(f_1,\ldots,f_{i-1})$ is 
injective, for all $i$. Then the Koszul complex 
$K^\bullet(A;f_1,\ldots,f_r)$ is acyclic in negative degrees. 

This is proved easily by induction on $r$, using (1.3), cf. [EGA III],  
Prop. 1.1.4. 

Obviously, 
$$
H^0(K^\bullet(A;f_1,\ldots,f_r))=A/(f_1,\ldots,f_r)
\eqno{(1.5)}
$$

{\bf 1.3.} From now on, until mentioned otherwise, the ring $A$ 
will be equal to a polynomial algebra $k[x_1,\ldots,x_n]$, 
$k$ being a fixed ground commutative ring. We denote by 
$$
\CT(A)=Der_k(A)
$$
the Lie algebra of $k$-linear derivations of $A$. We identify $\CT(A)$ 
with $A^n$ using the basis $\dpar_1,\ldots,\dpar_n$, where 
$\dpar_i:=\dpar/\dpar x_i$. 

We fix a function $f\in A$ such that $f(0)=0$. 
We denote by $M^\bullet(A;df)$ the Koszul complex
$$
M^\bullet(A;df)=K^\bullet(A;\dpar_1f,\ldots,\dpar_nf)
\eqno{(1.6)}
$$
Thus, 
$$
M^\bullet(A;df)=\Lambda^\bullet \CT(A)
\eqno{(1.7)}
$$
as a graded commutative algebra, with $\CT(A)$ having the degree $-1$. The  $(-1)-$component of the differential $\CT(A)\lra A$ sends 
$\tau\in\CT(A)$ to $\tau(f)$, cf. the remark at the end of 1.1. 

{\bf 1.4. Lemma.} {\it The following conditions are equivalent}: 

(i) $\dim(A/(\dpar_1f,\ldots,\dpar_nf))<\infty$; 

(ii) {\it $M^\bullet(A;df)$ is acyclic in negative degrees;} 

(iii) {\it $(\dpar_1f,\ldots,\dpar_nf)$ is a regular sequence.} 

This is a standard fact of Commutative Algebra, and follows from 
[S], Ch. IV, B.2, Th\'eor\`eme 2, Ch. III, B.3. 

If these equivalent conditions are fulfilled, we say that $f$ 
{\it has an isolated singularity} at $0$. The only non-zero cohomology 
space,
$$
A/(\dpar_1f,\ldots,\dpar_nf)=H^0(M^\bullet(A;df))
\eqno{(1.8)}
$$
is called the {\bf Milnor ring} of $f$ at $0$, cf. [M], \S 7.  

Thus, the commutative dg algebra $M^\bullet(A;df)$ may be considered as a 
resolution of the Milnor ring.                                                                                                                                                                                                                                                                                                                                                                                                                                                                                                                                                                                                                                                     

{\bf 1.5. Definition.} {\it A {\bf Gerstenhaber algebra} is a commutative 
$\BZ$-graded algebra $G^\bullet$, equipped with a bracket
$$
[\ ,\ ]:\ G^\bullet\otimes G^\bullet\lra G^\bullet[1]
$$
of degree $1$, which makes $G^\bullet[-1]$ a graded Lie algebra. 

The multiplication and the bracket must satisfy the following 
compatibility,
$$
[a,b\cdot c]=[a,b]\cdot c+(-1)^{(|a|-1)|b|}b\cdot [a,c]
\eqno{(1.9)}
$$
for all homogeneous $a,b, c$.}

Note the following equivalent form of (1.9), 
$$
[a\cdot b,c]=a\cdot [b,c]+(-1)^{(|c|-1)|b|}[a,c]\cdot b
\eqno{(1.10)}
$$
The equivalence of (1.9) and (1.10) is proved using the symmetry 
relations
$$
a\cdot b=(-1)^{|a||b|}b\cdot a
\eqno{(1.11)}
$$
and
$$
[a,b]=-(-1)^{(|a|-1)(|b|-1)}[b,a]
\eqno{(1.12)}
$$

{\bf 1.6. Theorem.} {\it There exists a unique bracket of degree $1$ 
on the algebra of polyvector fields $\Lambda^\bullet\CT(A)$, where to $\CT(A)$ 
the degree $-1$ and to $\CA$ the degree $0$ is assigned, such that} 

(a) {\it for $a\in A,\ \tau\in\CT(A),\ \ [\tau,a]=\tau(a)$;} 

(b) {\it the bracket on $\CT(A)$ coincides with the usual 
bracket of vector fields};   

(c) {\it together with the wedge 
product, this algebra becomes a Gerstenhaber algebra.} 

It is called the {\bf Schouten bracket}, [Sch]. The Gerstenhaber algebra 
defined in this Theorem, will be called the {\bf Schouten algebra}, 
and denoted $\CT^\bullet(A)$.  

Now, let us identify $\Lambda^\bullet\CT(A)$ with $M^\bullet(A;df)$. 

{\bf 1.7. Theorem.} {\it Equipped with the Schouten bracket, 
the complex $M^\bullet(A;df)[-1]$ is a dg Lie algebra.} 

We have to prove that 
$$
d([a,b])=[da,b]+(-1)^{|a|-1}[a,db]
\eqno{(1.13)}
$$
for $a\in\Lambda^{|a|}$. 
Using (1.12), it is easy to see that (1.13) is true for $(a,b)$ if and 
only if it is for $(b,a)$. 

{\bf 1.8. Lemma.} {\it If} (1.13) {\it holds true for $(a,b)$ and 
$(a,c)$ then it holds true for $(a,bc)$.}

This Lemma is proved by a straightforward computation, which 
we omit. One uses (1.10) and 
$$
d(a\cdot b)=da\cdot b+(-1)^{|a|}a\cdot db
\eqno{(1.14)}
$$
$\btu$

Returning to the Theorem, every element in $\Lambda^\bullet\CT(A)$ is a sum of 
elements, which are  
either functions, or products of vector fields. Therefore, due to the Lemma, 
we have to prove (1.13) only in the case when $a,b$ are either functions 
or vector fields. The only non-trivial identity comes out 
when both are vector fields. We have 
$$
d([\tau_1,\tau_2])=[\tau_1,\tau_2](f), 
$$
$\tau_i\in\CT(A)$. 
On the other hand, 
$$        
[d\tau_1,\tau_2]+(-1)^{-1-1}[\tau_1,d\tau_2]=[\tau_1(f),\tau_2]+
[\tau_1,\tau_2(f)]=
$$
$$
=-\tau_2\tau_1(f)+\tau_1\tau_2(f)=[\tau_1,\tau_2](f)
$$
The Theorem is proved. $\btu$

{\bf 1.9.} Let us call a {\bf degenerate Batalin-Vilkovisky algebra} a 
commutative dg algebra $G^\bullet$, together with a bracket 
$G^\bullet\otimes G^\bullet\lra G^\bullet[1]$ of degree $1$ making 
$G^\bullet[-1]$ a dg Lie algebra, such that with the 
differential forgotten, $G^\bullet$ is a Gerstenhaber algebra.  

The previous considerations show that the complex 
$M^\bullet(A;df)$, together with the wedge product of polyvector 
fields and the Schouten bracket, is a degenerate BV algebra.

\bigskip

\bigskip
\centerline{\bf \S 2. Batalin-Vilkovisky algebras on Calabi-Yau manifolds}
\bigskip

{\bf 2.1. Definition.} {\it A {\bf Batalin-Vilkovisky algebra} is a 
Gerstenhaber algebra $B^\bullet$, together with a differential $d$ 
of degree $1$ on it,  
such that 
$$
(-1)^{|a|}[a,b]=d(a\cdot b)-da\cdot b-(-1)^{|a|}a\cdot db
\eqno{(2.1)}
$$
for all homogeneous $a, b$.}

It follows from (2.1) that 
$$
d([a,b])=[da,b]+(-1)^{|a|-1}[a,db],
\eqno{(2.2)}
$$
i.e. $B^\bullet[-1]$ is a dg Lie algebra.  

{\bf 2.2.} Let us return to the situation of the last Section. Consider 
the de Rham complex
$$
\Omega^\bullet(A):\ 0\lra A\buildrel{d_{DR}}\over{\lra} 
\Omega^1(A)\buildrel{d_{DR}}\over{\lra}\ldots
\buildrel{d_{DR}}\over{\lra}\Omega^n(A)\lra 0
\eqno{(2.3)}
$$
Let us fix the top differential form
$$
\omega=dx_1\wedge\ldots\wedge dx_n
\eqno{(2.4)}
$$
Using $\omega$, one gets the isomorphisms
$$
\omega_*:\ \Lambda^p\CT(A)\iso\Omega^{n-p}(A)
\eqno{(2.5)}
$$
defined by
$$
\omega_*(a)=a\omega
\eqno{(2.6)}
$$
for $a\in A=\Lambda^0\CT(A)$, and 
$$
\omega_*(\tau_1\wedge\ldots\wedge\tau_p)=i_{\tau_1}\cdot\ldots\cdot i_{\tau_p}
(\omega)
\eqno{(2.7)}
$$
for $p\geq 1$. Explicitely, 
$$
\omega_*(a\dpar_{i_1}\wedge\ldots\wedge\dpar_{i_p})=
(-1)^{(\sum_{j=1}^p\ i_j)-p}adx_1\wedge\ldots\wedge\hat{dx}_{i_1}\wedge
\ldots\wedge\hat{dx}_{i_p}\wedge\ldots\wedge dx_n
\eqno{(2.8)}
$$
for $i_1<\ldots<i_p$. 

Let us identify the shifted de Rham complex $\Omega^\bullet(A)[n]$ 
with $\Lambda^\bullet\CT(A)$, using the isomorphisms (2.5). We can   
transfer the de Rham differential $d_{DR}$ to $\Lambda^\bullet\CT(A)$. 
Here $\CT(A)$ will have the degree $-1$. 

{\bf 2.3. Theorem.} {\it Equipped with the differential coming from 
$-d_{DR}$, the wedge product and the Schouten bracket, the 
algebra of polyvector fields $\Lambda^\bullet\CT(A)$ becomes a 
BV algebra.} 

We have to check the relation (2.1). First note that if (2.1) holds 
true for a pair $(a,b)$ then it does for $(b,a)$. 

{\bf 2.4. Lemma.} {\it Suppose that} (2.1) {\it is true for pairs 
$(a,b)$, $(b,c)$ and $(ab,c)$. Then it holds true the pair $(a,bc)$.} 

Indeed, we have 
$$
(-1)^{|a|}[a,bc]=(-1)^{|a|}[a,b]\cdot c+(-1)^{|a||b|+|a|+|b|}b\cdot [a,c],
\eqno{(2.9)}
$$
by (1.9). On the other hand, 
$$
d(a(bc))-da\cdot bc-(-1)^{|a|}a\cdot d(bc)=
$$
$$
=(-1)^{|a|+|b|}[ab,c]+d(ab)\cdot c+(-1)^{|a|+|b|}ab\cdot d(c)-
$$
$$
-da\cdot bc-(-1)^{|a|}a\cdot\bigl((-1)^{|b|}[b,c]+db\cdot c+
(-1)^{|b|}b\cdot dc\bigr)
$$
(we have applied (2.1) for $(ab,c)$ and $(b,c)$)
$$
=(-1)^{|a|+|b|}[ab,c]+\bigl((-1)^{|a|}[a,b]+da\cdot b+(-1)^{|a|}a\cdot db\bigr) 
\cdot c -
$$
$$
-da\cdot bc-(-1)^{|a|+|b|}a\cdot [b,c]-(-1)^{|a|}a\cdot db\cdot c
$$
(we have applied (2.1) for $(a,b)$)
$$
=(-1)^{|a|+|b|}[ab,c]+(-1)^{|a|}[a,b]\cdot c-
(-1)^{|a|+|b|}a\cdot [b,c]
\eqno{(2.10)}
$$
The right hand sides of (2.9) and (2.10) are equal, by the relation (1.10). 
This proves the Lemma. $\btu$

In view of this Lemma, to prove the Theorem, it is enough to check (2.1) 
when $a$ is either a function, or a vector field. Let us do that.  

If both $a$ and $b$ are functions, (2.1) is trivial.  


Let us check that 
$$
[a,b_1\dpar_{i_1}\wedge\ldots\wedge b_p\dpar_{i_p}]=
$$
$$
=
d(ab_1\dpar_{i_1}\wedge b_2\dpar_{i_2}\wedge\ldots\wedge\dpar_{i_p})-
a\cdot d(b_1\dpar_{i_1}\wedge\ldots\wedge b_p\dpar_{i_p}),
\eqno{(2.11)} 
$$
where $a\in A$, $i_1<\ldots<i_j$. 
Using the definition, one checks that
$$
d(b_1\dpar_{i_1}\wedge\ldots\wedge b_p\dpar_{i_p})=
\sum_{j=1}^p\ (-1)^j\dpar_{i_j}(b_1\cdot\ldots\cdot b_p)\dpar_{i_1}\wedge
\ldots\wedge\hat{\dpar}_{i_j}\wedge\ldots\wedge\dpar_{i_p}
\eqno{(2.12)}
$$
On the other hand, the repeated application of (1.9) gives
$$
[a,b_1\dpar_{i_1}\wedge\ldots\wedge b_p\dpar_{i_p}]=
\sum_{j=1}^p\ (-1)^jb_j\dpar_{i_j}(a)b_1\dpar_{i_1}\wedge\ldots\wedge
\widehat{b_j\dpar}_{i_j}\wedge\ldots\wedge b_p\dpar_{i_p}
\eqno{(2.13)}
$$
These two formulas easily imply (2.11). This completes the check of 
(2.1) for $a\in A$. 

The check of (2.1) for $a\in\CT(A)$ is similar; we leave it to the reader. 
The Theorem is proven. $\btu$

{\bf 2.5.} Now let us turn the function $f$ on. First of all, note that 
when we identify $\Omega^\bullet(A)$ with $\Lambda^\bullet\CT(A)$ using 
$\omega_*$ as in (2.2), the multiplication by $df$ becomes 
the Koszul differential in $M^\bullet(A;df)$, (1.6). 

{\bf 2.6. Theorem.} {\it The Schouten algebra $\CT^\bullet(A)$} (cf. 1.6) 
{\it equipped with the differential coming from $-d_{DR}+df$, 
is a Batalin-Vilkovisky algebra.} 

{\bf 2.7. Lemma.} {\it Let $S^\bullet$ be a Gerstenhaber algebra. Let 
$d_0, d_1$ be two anticommuting differentials on $S^\bullet$, such that 
$(S^\bullet,d_0)$ is a degenerate BV-algebra} (cf. 1.9), 
{\it and $(S^\bullet,d_1)$ is 
a BV-algebra. Then $(S^\bullet,d_0+d_1)$ is a BV-algebra.} 

The proof is obvious. The anticommutation of $d_0$ and $d_1$ is needed 
for the identity $(d_0+d_1)^2=0$. 
$\btu$

Applying the Lemma to the Schouten algebra $\CT^\bullet(A)$, to $d_0$ 
induced from $df$, and $d_1$ induced from $-d_{DR}$, we get the Theorem, 
in view of 2.5 and Theorems 1.7 and 2.3. $\btu$   

The BV-algebra defined in Theorem 2.6 will be denoted $\CM(A,\omega;df)$. 

{\bf 2.8.} Let us make a base change $k\lra k[\lambda]$ where 
$\lambda$ is an independent variable. Using 2.6, we get a BV-algebra over 
$k[\lambda]$, $\CM(A[\lambda],\omega;\lambda\cdot df)$. 
The degenerate BV-algebra $M^\bullet(A;df)$ may be regarded as a "limit" of 
$\CM(A,\omega;\lambda\cdot df)$ at $\lambda\lra\infty$.

{\bf 2.9.} We can apply the same construction by taking as an algebra $A$ 
a localization of the polynomial algebra $k[x_1,\ldots,x_n]$, and 
as a form $\omega$ an arbitrary generator
$$
\omega=c(x)dx_1\wedge\ldots\wedge dx_n,
\eqno{(2.14)}
$$
of the $A$-module $\Omega^n(A)$. Here $c(x)$ is a unit  
in $A$. 

Theorem 2.3 remains true, with the same proof. We should 
apply the following modifications of the formulas (2.11) and (2.13). 
$$
d(a\dpar_i)=-\dpar_i(a)-a\cdot \dpar_i\log c,
\eqno{(2.15)}
$$
and
$$
d(a\dpar_i\wedge b\dpar_j)=-\bigl(\dpar_i(ab)+ab\cdot\dpar_i\log c\bigr)\dpar_j+
\bigl(\dpar_j(ab)+ab\cdot\dpar_j\log c)\dpar_i
\eqno{(2.16)}
$$
Here
$$
\dpar_i\log c=\dpar_ic\cdot c^{-1}
\eqno{(2.17)}
$$
Theorem 2.6 remains true, with the same proof. 

{\bf 2.10.} Let $A$ be an arbitrary commutative $k$-algebra, 
smooth over $k$. Assume that $A$ is \'etale over a polynomial algebra 
$B=k[x_1,\ldots,x_n]$. By definition, cf. [SGA 1], Expos\'e II, 
D\'ef. 1.1, this holds true 
for any smooth $k$-algebra, locally over $Spec(A)$.
 
The derivations $\dpar/\dpar_{x_i}\in \CT(B)$ extend uniquely to derivations 
$\dpar_i\in \CT(A)$, which form a basis of $\CT(A)$ as an $A$-module, 
and anticommute. Their duals $dx_i$ form a basis of $\Omega^n(A)$.  

Choose a generator
$$
\omega=cdx_1\wedge\ldots\wedge dx_n
\eqno{(2.17)}
$$
of $\Omega^n(A)$. It induces the isomorphisms
$$
\omega_*:\ \Lambda^\bullet\CT(A)\iso\Omega^{n-\bullet}(A)
\eqno{(2.18)}
$$
Theorems 2.3 and 2.6 remain true, with the same proof. This follows from 2.9. 

{\bf 2.11.} Now let us pass to the global situation. Suppose we are given 
a scheme $X$, smooth of relative dimension $n$ over $k$. Assume that the 
canonical bundle $\Omega^n_X$ is trivial. 

Let us choose a trivialization
$$
\omega:\ \CO_X\iso\Omega^n_X
\eqno{(2.19)}
$$
It induces the isomorphisms of graded $\CO_X$-modules
$$
\omega_*:\ \Lambda^\bullet\CT_X\iso\Omega^{n-\bullet}_X
\eqno{(2.20)}
$$

{\bf 2.12. Theorem.} {\it Given a closed one-form 
$\phi\in\Gamma(X;\Omega^1_X)$, consider the Schouten algebra 
$\Lambda^\bullet\CT_X$, equipped with the differential obtained, using the 
identification} (2.18), {\it from the differential 
$$
-d_{DR}+\phi
\eqno{(2.21)}
$$
on $\Omega^\bullet_X$.  
This way we get a sheaf of Batalin-Vilkovisky algebras on $X$.}

Indeed, the statement is local, and we can assume we are in the situation 
2.10. Now, consider the situation of Theorem 1.7, with $df$ replaced 
by an {\it arbitrary} one-form $\phi\in\Omega^1(A)$. In the corresponding 
complex $M^\bullet(A;\phi)$, the differential is equal to the  
multiplication by $\phi$, if we use the identification (2.18). 

In particular the differential $\CT(A)\lra A$ is given by
$$
\tau\mapsto\langle\tau,\phi\rangle
\eqno{(2.22)}
$$
The conclusion of Theorem 1.7 remains true in these more general assumptions. 
Indeed, we can use the same Lemma 1.8, and replace the last part of the proof 
of 1.7 by the Lemma below. 

{\bf 2.13. Lemma.} {\it For all $\tau_1, \tau_2\in\CT(A)$, and 
$\phi\in\Omega^1(A)$, we have 
$$
\langle [\tau_1,\tau_2],\phi\rangle=\tau_1(\langle\tau_2,\phi\rangle)-
\tau_2(\langle\tau_1,\phi\rangle)
\eqno{(2.23)}
$$}

In fact, we can assume that $\tau_1=a\dpar_i, \tau_2=b\dpar_j$ and 
$\phi=dx_k$. Then the left hand side of (2.23) is equal to 
$$
\langle a\dpar_i(b)\dpar_j-b\dpar_j(a)\dpar_i,dx_k\rangle=
a\dpar_i(b)\delta_{jk}-b\dpar_j(a)\delta_{ik},
$$
which is equal to the right hand side of (2.23). $\btu$

Using this generalization of Theorem 1.7, we get a generalization of 
Theorem 2.6, with $df$ replaced by an arbitrary {\it closed} one 
form $\phi$, and $A$ as in 2.10. (The closedness of $\phi$ is necessary 
and sufficient condition for the equality $(-d_{DR}+\phi)^2=0$.) 
This generalization of 2.6 implies, in turn, Theorem 2.12. $\btu$     

The sheaf of BV algebras defined in this Theorem, will be denoted 
$\CM^\bullet(X,\omega;\phi)$.

\bigskip

\bigskip
\centerline{\bf \S 3. Recollections on $D$-modules}
\bigskip

For more details about the generalities below, cf. [BD1], Chapter 2, \S 1,  
and [B].  

{\bf 3.1.} From now on, we assume that our ground ring $k$ is a field of 
characteristic $0$. Let $X$ be a scheme, smooth over $k$, of dimension $n$. 
We will consider the categories $Mod(\CD_X)^\ell$ and $Mod(\CD_X)^r$ of left 
(resp., right) $\CD_X$-modules ($:=$ sheaves of left (resp., right) 
$\CD_X$-modules, quasicoherent over $\CO_X$).

The sheaf of rings $\CD_X$ is generated by two $\CO_X$-modules, $\CO_X$ and 
$\CT_X$, and relations
$$
\tau\cdot a-a\cdot\tau=\tau(a)
\eqno{(3.1)}
$$
and
$$
\tau\cdot\tau'-\tau'\cdot\tau=[\tau,\tau'],
\eqno{(3.2)}
$$
for $\tau, \tau'\in\CT_X,\ a\in\CO_X$. 

Therefore, a left $\CD_X$ module is a quasicoherent $\CO_X$-module $M$, 
equipped with the left multiplication by $\CT_X$, such that 
$$
a(\tau m)=(a\tau)\cdot m
\eqno{(3.3a)}
$$
$$
\tau (am)=a(\tau m)+\tau(a)\cdot m
\eqno{(3.3b)}
$$
$$
\tau(\tau'm)-\tau'(\tau m)=[\tau,\tau']\cdot m,
\eqno{(3.3c)}
$$
for $m\in M$. 

Similarly, a right $\CD_X$-module is a quasicoherent $\CO_X$-module $N$, 
equipped with the right multiplication by $\CT_X$, such that
$$
(na)\tau=n\cdot(a\tau)
\eqno{(3.4a)}
$$
$$
(n\tau)a=n(a\tau)+n\cdot\tau(a)
\eqno{(3.4b)}
$$
$$
(n\tau)\tau'-(n\tau')\tau=n\cdot [\tau,\tau'],
\eqno{(3.4c)}
$$
for $n\in N$. Here we agree that 
$$
n\cdot a=a\cdot n
\eqno{(3.5)}
$$ 
    
{\bf 3.2.} For $M, M'\in Mod(\CD_X)^\ell$, the tensor product 
$M\otimes_{\CO_X}N$ admits a canonical structure of a left $\CD_X$-module, 
defined by the Leibnitz rule
$$
\tau\cdot(m\otimes m')=\tau m\otimes m'+m\otimes\tau m',
\eqno{(3.6)}
$$
for $\tau\in\CT_X$. This way we get a tensor structure on $Mod(\CD_X)^\ell$ 
(it is denoted by $\otimes^!$ in [BB], cf. {\it op. cit.}, 2.2.1).

The sheaf $Hom_{\CO_X}(M,M')$ is also canonically a left $\CD_X$-module, 
with the $\CT_X$-action given by
$$
(\tau f)(m)=\tau(f(m))-f(\tau m),
\eqno{(3.7)}
$$
for $f\in Hom_{\CO_X}(M,M')$. The rule (3.7) is invented in such a way that 
the canonical morphism, 
$$
Hom_{\CO_X}(M,M')\otimes M\lra M'
\eqno{(3.8)}
$$
would be a morphism of $\CD_X$-modules.

{\bf 3.3.} For $M\in Mod(\CD_X)^\ell,\ N\in Mod(\CD_X)^r$, the tensor 
product $M\otimes N$ is canonically a right $\CD_X$-module, with the 
multiplication by $\CT_X$ defined by
$$
(m\otimes n)\cdot\tau=m\otimes n\tau-\tau m\otimes n
\eqno{(3.9)}
$$
For $N'\in Mod(\CD_X)^r$, the $\CO_X$-module $Hom_{\CO_X}(N,N')$ is 
canonically a left $\CD_X$-module, with the multiplication by $\CT_X$ 
given by 
$$
(\tau f)(n)=f(n\tau)-f(n)\cdot \tau
\eqno{(3.10)}
$$
Again, this rule is written in such a way that the canonical 
adjunction morphism is a morphism of right $\CD_X$-modules. 

{\bf 3.4.} The structure sheaf $\CO_X$ has an obvious structure of a 
left $\CD_X$-module. 

The canonical bundle $\Omega^n_X$ is canonically a right $\CD_X$-module. 
The multiplication by $\CT_X$ is given by 
$$
\omega\cdot\tau=-Lie_{\tau}(\omega),
\eqno{(3.11)}
$$
for $\omega\in\Omega^n_X$. 

We have the {\bf Cartan formula}
$$
i_{\tau}\circ d+d\circ i_{\tau}=Lie_{\tau}
\eqno{(3.12)}
$$
Consequently, 
$$
\omega\cdot\tau=-d(i_{\tau}\omega)
\eqno{(3.13)}
$$
Let us check the relations (3.4a-c). We have 
$$
\omega\cdot (a\tau)=-d i_{a\tau}(\omega)=-di_{\tau}(a\omega)=
(\omega a)\cdot\tau, 
$$
which proves (3.4a). We have 
$$
\omega\cdot (a\tau)=-Lie_{\tau}(a\omega)=-\tau(a)\omega-a Lie_{\tau(\omega)}, 
$$
which proves (3.4b). The relation (3.4c) is obvious. 

Using 3.3, we get the functors 
$$
Mod(\CD_X)^\ell\lra Mod(\CD_X)^r,\ M\mapsto M\otimes\Omega^n_X, 
\eqno{(3.14)}
$$
and 
$$
Mod(\CD_X)^r\lra Mod(\CD_X)^\ell,\ N\mapsto Hom(\Omega^n_X,N)
\eqno{(3.15)}
$$
which are mutually inverse equivalences. 

{\bf 3.5. De Rham complex.} We have the De Rham complex of our scheme $X$, 
$$
\Omega^\bullet_X:\ 0\lra\CO_X\lra\Omega^1_X\lra\ldots\lra\Omega^n_X\lra 0
\eqno{(3.16)}
$$
We can regard $\Omega^\bullet_X$ as a twisted version of the standard 
cochain complex of the Lie algebra of vector fields $\CT_X$ with coefficients 
in the module of functions $\CO_X$. Namely, let us identify
$$
\Omega^i_X=Hom_{\CO_X}(\Lambda^i\CT_X,\CO_X)
\eqno{(3.17)}
$$
Thus, we get the identification
$$
\Omega^\bullet_X=Hom_{\CO_X}(\Lambda^\bullet\CT_X,\CO_X),
\eqno{(3.18)}
$$
with the differential defined by
$$
d\omega(\tau_1\wedge\ldots\wedge\tau_p)=\sum_{i=1}^p\ (-1)^{i+1}\tau_i
(\omega(\tau_1\wedge\ldots\wedge\hat{\tau}_i\wedge\ldots\wedge\tau_p))+
$$
$$
+\sum_{1\leq i<j\leq p}\ (-1)^{i+j+1}\ \omega([\tau_i,\tau_j]\wedge
\tau_1\wedge\ldots\wedge\hat{\tau}_i\wedge\ldots\wedge
\hat{\tau}_j\wedge\ldots\wedge\tau_p),
\eqno{(3.19)}
$$
for $\omega\in Hom(\Lambda^p\CT_X,\CO_X)$. 
cf. [KN], Chapter I, Prop. 3.11.  

Twisted, since the bracket in $\CT_X$ is not $\CO_X$-linear. The differentials 
are not $\CO_X$-linear, but are differential operators of the 
first order. 
 
In the same sense, $\CD_X$ is the twisted enveloping algebra of $\CT_X$, 
$$
\CD_X=U_{\CO_X}\CT_X
\eqno{(3.20)}
$$   
More generally, for a left $\CD_X$-module $M\in Mod(\CD_X)^\ell$, we have 
the De Rham complex
$$
DR^\ell(M):\ 
0\lra M\lra \Omega^1_X\otimes M\lra\ldots\lra\Omega^n_X\otimes M\lra 0
\eqno{(3.21)}
$$
Using the identifications (3.17), we can identify this complex with 
the twisted cochain complex of the Lie algebra $\CT_X$ with coefficients 
in $M$, 
$$
DR^\ell(M)=Hom_{\CO_X}(\Lambda^\bullet\CT_X,M),
\eqno{(3.22)}
$$
the differential being defined by the same formula (3.19). 

For a right $\CD_X$-module $N\in Mod(\CD_X)^r$, we have the  
De Rham complex 
$$
DR(N)^r:\ 
0\lra N\otimes\Lambda^n\CT_X\lra\ldots\lra N\otimes\CT_X\lra N\lra 0
\eqno{(3.23)}
$$
We consider it as living in degrees from $-n$ to $0$.  
It may be regarded as the twisted version of the standard chain complex 
of the Lie algebra $\CT_X$ with coefficients in $N$, 
$$
DR(N)^r=N\otimes_{\CO_X}\Lambda^\bullet\CT_X
\eqno{(3.24)}
$$ 
The differential acts as follows
$$
d(n\otimes\tau_1\wedge\ldots\wedge\tau_p)=\sum_{i=1}^p\ (-1)^{i-1}
n\tau_i\otimes\tau_1\wedge\ldots\wedge\hat{\tau}_i\wedge\ldots\wedge 
\tau_p+
$$
$$
+\sum_{1\leq i<j\leq p}\ (-1)^{i+j}n\otimes [\tau_i,\tau_j]\wedge
\ldots\wedge\hat{\tau}_i\wedge\ldots\wedge\hat{\tau}_j\wedge\ldots
\wedge\tau_p
\eqno{(3.25)}
$$
Let us check that (3.25) is a well defined map. For example, we have 
to check that 
$$
d(n\otimes a\tau_1\wedge\tau_2)=d(n\otimes\tau_1\wedge a\tau_2), 
\eqno{(3.26)}
$$
for $a\in\CO_X$. We have
$$
d(n\otimes a\tau_1\wedge\tau_2)=na\tau_1\otimes\tau_2-
n\tau_2\otimes a\tau_1-n\otimes [a\tau_1,\tau_2]=
$$
$$
=n\tau_1 a\otimes\tau_2-n\tau_1(a)\otimes\tau_2-n\tau_2\otimes a\tau_1-
n\otimes a[\tau_1,\tau_2]+n\otimes\tau_2(a)\tau_1
$$
We have used the formula
$$
[a\tau_1,\tau_2]=a[\tau_1,\tau_2]-\tau_2(a)\tau_1
\eqno{(3.27)}
$$
which is a particular case of (1.10). On the other hand, 
$$
d(n\otimes\tau_1\wedge a\tau_2)=n\tau_1\otimes a\tau_2-na\tau_2\otimes\tau_1-
n\otimes [\tau_1,a\tau_2]=
$$
$$
=n\tau_1\otimes a\tau_2-n\tau_2 a\otimes\tau_1+n\tau_2(a)\otimes\tau_1-
n\otimes a[\tau_1,\tau_2]-n\otimes\tau_1(a)\tau_2
$$
which is equal to the left hand side of (3.26). In the same manner, 
the analog of (3.26) with an arbitrary number of $\tau$'s is checked. 
This proves that the map (3.25) is well defined. 

The complexes $DR^\ell$ and $DR^r$ are mapped to each other under the 
canonical equivalence between left and right $\CD_X$-modules, 
(3.14), (3.15). 

\bigskip

\bigskip
{\bf \S 4. Connections and Batalin-Vilkovisky structures}
\bigskip

{\bf 4.1.} We keep the assumptions of the last Section. Let us call 
a {\bf Calabi-Yau data} on $X$ an integrable connection 
$$
\nabla:\ \Omega^n_X\lra\Omega^1_X\otimes\Omega^n_X
\eqno{(4.1)}
$$
on the canonical bundle $\Omega^n_X$. 

In other words, a CY data is a structure of a left $\CD_X$-module 
on $\Omega^n_X$. By 3.4, this is the same as a structure of a 
{\it right} $\CD_X$-module on the structure sheaf 
$\CO_X=Hom_{\CO_X}(\Omega^n_X,\Omega^n_X)$, i.e. 
a map 
$$
\nabla^r:\ \CT_X=\CO_X\otimes_{\CO_X}\CT_X\lra\CO_X
\eqno{(4.2)}
$$
such that 
$$
\nabla^r(a\tau)=a\nabla^r(\tau)-\tau(a),
\eqno{(4.3)}
$$
for $a\in\CO_X,\ \tau\in\CT_X$. 

{\bf 4.2.} Let us call a {\bf Batalin-Vilkovisky data} on $X$ a differential 
on the Schouten algebra $\Lambda^\bullet\CT_X$, making it a BV algebra. 

{\bf 4.3. Theorem.} {\it Given a CY data on $X$, the De Rham complex 
$DR^r(\CO_X)$ of  
the corresponding right $\CD_X$-module $\CO_X$ is a 
BV data on $X$. 

Conversely, given a BV data on $X$, the component $\CT_X\lra\CO_X$ 
of the corresponding differential on $\Lambda^\bullet\CT_X$ gives a map} 
(4.2) {\it which defines a CY data on $X$. 

These two procedures establish the mutually inverse bijections between 
the sets of CY data and BV data on $X$.}

Suppose we are given a CY data. We have to show the $DR^r(\CO_X)$ is a 
BV algebra. This is a generalization of Theorem 2.12. We have to check 
the identity (2.1). We will use the formula
$$
d(\tau_1\wedge\ldots\wedge\tau_p)=\sum_{j=1}^p\ (-1)^{j-1}
\nabla^r(\tau_j)\tau_1\wedge\ldots\wedge\hat{\tau}_j\wedge\ldots\wedge
\tau_p+
$$
$$
+\sum_{i<j}\ (-1)^{i+j}[\tau_i,\tau_j]\wedge\tau_1\wedge\ldots\wedge
\hat{\tau}_i\wedge\ldots\wedge\hat{\tau}_j\wedge\ldots\wedge\tau_p
\eqno{(4.4)}
$$
which is a version of (3.25). 

Application of Lemma 2.4 shows that it is 
enough to check (2.1) when $a$ is either a function or a vector field. 
If both $a, b$ are functions, (2.1) is trivial. 

Let $a\in\CO_X$. We have
$$
[a,\tau_1\wedge\ldots\wedge\tau_p]=\sum_{j=1}^p\ (-1)^p\tau_j(a)
\tau_1\wedge\ldots\wedge\hat{\tau_j}\wedge\ldots\wedge\tau_p
$$
On the other hand, it follows from (4.4) that  
$$
d(a\tau_1\wedge\ldots\wedge\tau_p)=\nabla^r(a\tau_1)\wedge\tau_2
\wedge\ldots\wedge\tau_p+
\sum_{j=2}^p\ (-1)^{j-1}
\nabla^r(\tau_j)a\tau_1\wedge\ldots\wedge\hat{\tau}_j\wedge\ldots\wedge
\tau_p+
$$
$$
+\sum_{1<j}\ (-1)^{1+j}[a\tau_1,\tau_j]\wedge\tau_2\wedge\ldots\wedge
\hat{\tau}_j\wedge\ldots\wedge\tau_p+
\sum_{1<i<j}\ (-1)^{i+j}a[\tau_i,\tau_j]\wedge\tau_1\wedge\ldots\wedge
\hat{\tau}_i\wedge\ldots\wedge\hat{\tau}_j\wedge\ldots\wedge\tau_p=
$$
$$
=\bigl(a\nabla^r(\tau_1)-\tau_1(a)\bigr)\wedge\tau_2\wedge\ldots\wedge\tau_p+
\sum_{j=2}^p\ (-1)^{j-1}
a\nabla^r(\tau_j)\tau_1\wedge\ldots\wedge\hat{\tau}_j\wedge\ldots\wedge
\tau_p+
$$
$$
+\sum_{1<j}\ (-1)^{1+j}\bigl(a[\tau_1,\tau_j]-\tau_j(a)\bigr)\wedge\tau_2\wedge\ldots\wedge
\hat{\tau}_j\wedge\ldots\wedge\tau_p+
$$
$$
+\sum_{1<i<j}\ (-1)^{i+j}a[\tau_i,\tau_j]\wedge\tau_1\wedge\ldots\wedge
\hat{\tau}_i\wedge\ldots\wedge\hat{\tau}_j\wedge\ldots\wedge\tau_p
$$
Subtracting from this $a\cdot d(\tau_1\wedge\ldots\wedge\tau_p)$, we get (2.1) 
in this case.

It remains to check the following particular case of (2.1), 
$$
-[\tau_0,\tau_1\wedge\ldots\wedge\tau_p]=d(\tau_0\wedge\ldots\wedge\tau_p)-
d\tau_0\wedge\tau_1\wedge\ldots\wedge\tau_p+
\tau_0\wedge d(\tau_1\wedge\ldots\wedge\tau_p)
\eqno{(4.5)}
$$
The left hand side of (4.5) is equal to
$$
\sum_{j=1}^p\ (-1)^j[\tau_0,\tau_j]\wedge\tau_1\wedge\ldots\wedge\hat{\tau}_j
\wedge\ldots\wedge\tau_p
$$
The right hand side is equal to
$$
\sum_{j=0}^p\ (-1)^j\nabla^r(\tau_j)\wedge\tau_0\wedge\ldots\wedge\hat{\tau}_j
\wedge\ldots\wedge\tau_p+
\sum_{1\leq j}\ (-1)^j[\tau_0,\tau_j]\wedge\tau_0\wedge\ldots\wedge\hat{\tau}_j
\wedge\ldots\wedge\tau_p+
$$
$$
+\sum_{1\leq i<j}\ (-1)^{i+j}[\tau_i,\tau_j]\wedge\tau_0\wedge\ldots\wedge
\hat{\tau}_i\wedge\ldots\wedge\hat{\tau}_j\wedge\ldots\wedge\tau_p-
\nabla^r(\tau_0)\wedge\tau_1\wedge\ldots\wedge\tau_p+
$$
$$
+\sum_{j=1}^p\ (-1)^{j-1}\nabla^r(\tau_j)\tau_0\wedge\ldots\wedge
\hat{\tau}_j\wedge\ldots\wedge\tau_p+
\sum_{1\leq i<j}\ (-1)^{i+j}\tau_0\wedge [\tau_i,\tau_j]\wedge\tau_1
\wedge\ldots\wedge\hat{\tau}_i\wedge\ldots\wedge\hat{\tau}_j\wedge\ldots\wedge
\tau_p
$$
which is equal to the left hand side. This proves (4.5) and completes 
the construction of the arrow $\alpha:\ $(CY data) $\lra$ (BV data). 

Conversely, given a BV data, a particular case of (2.1) is
$$
[a,\tau]=\nabla^r(a\tau)-a\nabla^r(\tau)
$$
which is nothing but (4.3). This gives the arrow
$\beta:\ $ (BV data) $\lra$ (CY data). It is clear that the composition 
$\beta\alpha$ is the identity. On the other hand, given a BV data, 
the differential in $\Lambda^\bullet\CT_X$ is uniquely determined 
by its last component, and (2.1). This implies that $\alpha\beta$ is the 
identity. 

The Theorem is proven. $\btu$

{\bf 4.4.} Let $\CA$ be a commutative $\CD_X$-algebra, i.e. a commutative 
(implied: associative, unital)  
algebra in the category of $Mod(\CD_X)^\ell$ , with the $\otimes^!$ 
tensor structure, cf. 3.2. On the sheaf of $\CA$-modules
$$
\CT_{\CA}:=\CA\otimes_{\CO}\CT_X,
\eqno{(4.6)}
$$
one defines a Lie algebra structure by
$$
[a\tau,b\tau']=a\tau(b)\tau'-b\tau'(a)\tau+ab [\tau,\tau'], 
\eqno{(4.7)}
$$
Here by $a\tau$ we denoted $a\otimes\tau$, for $a\in\CA,\ \tau\in\CT_X$, 
and by $\tau(a)$ we denoted $\tau\cdot a\in\CA$, defined by  the 
left $\CD_X$-module structure on $\CA$. 

The given $\CO_X$-linear map $\CT_X\lra Der(\CA)$ is extended by linearity 
to the $\CA$-linear map
$$
\CT_\CA\lra Der(\CA)
$$
which defines on $\CA$ a structure of a module over the Lie algebra 
$\CT_\CA$. The action of $\CT_\CA$ on $\CA$ will also be denoted by 
$\tau(a)$, $a\in\CA,\ \tau\in \CA$. We have 
$$
[\tau,a\tau']=a[\tau,\tau']+\tau(a)\tau'
\eqno{(4.8)}
$$
for all $\tau, \tau'\in\CT_\CA,\ a\in\CA$. Thus, $\CT_\CA$ is a 
{\it Lie algebroid over $\CA$}, in the terminology of [BD2], 3.5.10.   

Let $\Lambda^\bullet(\CT_{\CA})$ denote the exterior algebra of 
$\CT_{\CA}$ {\it over $\CA$}. We regard it as a graded commutative 
algebra, with $\CA$ (resp. $\CT_{\CA}$) having degree $0$ (resp. $-1$).   

We have the following generalization of Theorem 1.6. 

{\bf 4.5. Theorem.} {\it There exists a unique Lie bracket of degree $1$ 
on $\Lambda^\bullet(\CT_{\CA})$ such that}

(a) {\it for $a\in\CA,\ \tau\in\CT_{\CA},\ [\tau,a]=\tau(a)$}; 

(b) {\it the bracket on $\CT_{\CA}$ is given by} (4.7); 

(c) {\it equipped with this bracket, $\Lambda^\bullet(\CT_{\CA})$ 
becomes a Gerstenhaber algebra.}

{\bf 4.6.} Let $\CD_\CA$ denote the sheaf of associative algebras 
generated by the algebra $\CA$ and the sheaf $\CT_X$, subject 
to relations 

(a) $\tau\cdot a-a\cdot\tau=\tau(a)$; 

(b) $\tau\cdot\tau'-\tau'\cdot\tau=[\tau,\tau']$.

We have a canonical filtration $\{\CD^{\leq i}_{\CA}\}$ on $\CD_\CA$, 
with the graded algebra
$$
gr^\bullet\CD_\CA=S^\bullet(\CT_\CA),
\eqno{(4.9)}
$$
the symmetric algebra over $\CA$. 

One has a canonical isomorphism of left $\CA$-modules
$$
\CD_\CA=\CA\otimes_\CO\CD_X
\eqno{(4.10)}
$$
We have canonical embeddings of 
$\CO_X$-algebras
$$
\CD_X\hra\CD_\CA;\ \ \CA\hra\CD_\CA
$$
(We have just described the {\it enveloping algebra} construction from 
[BB], 1.2.5.) 

{\bf 4.7.} The algebra $\CA$ has an obvious structure of a left 
$\CD_\CA$-module. Assume we are given a structure of a right $\CD_\CA$-module 
on $\CA$ which induces an obvious $\CA$-module structure. 
After restriction, $\CA$ becomes a right $\CD_X$-module, such that 
$$
(a\cdot b)\cdot\tau=a\cdot(b\cdot\tau)-b\cdot\tau(a),
\eqno{(4.11)}
$$
for all $a,b\in\CA,\ \tau\in\CT_X$. 
Define a map $\nabla^r:\ \CT_\CA\lra\CA$ by
$$
\nabla^r(a\tau)=a\cdot\tau
\eqno{(4.12)}
$$
We have
$$
\nabla^r(a\nu)=a\nabla^r(\nu)-\nu(a)
\eqno{(4.13)}
$$
for all $a\in\CA,\ \nu\in\CT_\CA$. 

{\bf 4.8. Lemma.} {\it The previous definitions establish bijections 
between the data} (i) --- (iv) {\it below.} 

(i) {\it The structures of a right $\CD_\CA$-module on $\CA$, inducing the 
obvious $\CA$-module structure.}

(ii) {\it The structures of a right $\CD_X$-module on $\CA$ inducing the given 
$\CO_X$-module structure, such that} (4.11) {\it is satisfied.}

(iii) {\it The maps}
$$
\nabla^r:\ \CT_X\lra\CA
$$
{\it such that} (4.13) {\it is satisfied, for all $a\in\CO_X,\ \nu\in\CT_X$.}

(iv) {\it The maps
$$
\nabla^r:\ \CT_\CA\lra\CA
$$
such that} (4.13) {\it is satisfied.} 

A structure described in this Lemma will be called a 
{\bf CY$_{\CA}$-structure}. 

{\bf 4.9. Theorem.} {\it Given a CY$_\CA$-structure, there exists a unique 
differential of degree $1$ on the Gerstenhaber algebra 
$\Lambda^\bullet\CT_\CA$,  
whose first component coincides with} (4.13) {\it 
and making $\Lambda^\bullet\CT_\CA$ a Batalin-Vilkovisky algebra.}

Let us call a differential on $\Lambda^\bullet\CT_\CA$ making it a BV 
algebra a {\bf BV$_\CA$-structure}. Given a BV$_\CA$-structure, the 
first component of the differential is a map 
$\CT_\CA\lra\CA$ satisfying (4.13). 

{\bf 4.10. Theorem.} {\it The above constructions establish the two 
inverse bijections between the set of CY$_\CA$-data and the set of 
BV$_\CA$-data.} 

This is a generalization of Theorem 4.3. The proof is the same. On just 
notes that in the proof of 4.3, we have used only the identity (4.13).  

{\bf 4.11.} What is a Gerstenhaber algebra (over $X$) living in degrees 
$-1, 0$? It is the same as a sheaf of commutative algebras $\CA^0$, 
a sheaf of Lie algebras $\CA^{-1}$, which is also 
an $\CA^0$-module, and acts on $\CA^0$ by derivations. The corresponding 
map $\CA^{-1}\lra Der(\CA^0)$ should be $\CA^0$-linear, 
and the identity (4.8) should hold
for all $a\in\CA^0,\ \tau\in\CA^{-1}$. In other words, 
$\CA^{-1}$ is a Lie algebroid over $\CA^0$. 

Let $Gerst$ denote the category of Gerstenhaber algebras over $X$, and let  
$$
Gerst^{[-1,0]}\subset Gerst
\eqno{(4.14)}
$$ 
be the full subcategory consisting 
of Gerstenhaber algebras living in degrees $-1, 0$. We claim that the 
obvious trucation functor
$$
t:\ Gerst\lra Gerst^{[-1,0]}
\eqno{(4.15)}
$$ 
admits a left adjoint 
$$
S:\ Gerst^{[-1,0]}\lra Gerst
\eqno{(4.16)}
$$ 
In fact, given $\CA\in Gerst^{[-1,0]}$, set 
$$
S(\CA)=S^\bullet_{\CA^0}(\CA^{-1}[1]),
\eqno{(4.17)}
$$
the symmetric algebra over $\CA^0$. 
In other words, as a graded algebra 
$$
S(A)=\Lambda^\bullet_{\CA^0}(\CA^{-1}),
\eqno{(4.18)}
$$
with $\CA^{-1}$ having degree $-1$. 

{\bf 4.12. Theorem.} {\it There exists a unique Lie bracket of degree $1$ 
on $S(\CA)$, making it a Gerstenhaber algebra, such that the 
obvious embedding $\CA\hra S(\CA)$ is a morphism of graded Lie algebras.} 

This is a version of Theorems 4.5 and 1.6. This bracket is called 
the {\it Schouten bracket}. This way we get a structure of a Gerstenhaber 
algebra on $S(\CA)$, called the {\bf Schouten algebra of} $\CA$.  

The functor $S$ is the left adjoint to (4.15).  
Obviously, the composition $t\circ S$ is the identity.  

{\bf 4.13.} Given $\CA\in Gerst^{[-1,0]}$, assume we have a BV structure 
on $S(\CA)$, i.e. a differential of degree $1$ on this algebra, making 
it a BV algebra. Its first component is a mapping 
$$
\nabla^r:\ \CA^{-1}\lra\CA^0
\eqno{(4.19)}
$$
such that
$$
\nabla^r(a\tau)=a\nabla^r(\tau)-\tau(a)
\eqno{(4.20)}
$$
for all $a\in\CA^0,\ \tau\in\CA^{-1}$. 

Conversely, given a mapping (4.19) satisfying (4.20), there exists a unique 
extension of this mapping to a differential $d$ on $S(\CA)$, making it 
a BV algebra. In fact, the uniqueness of $d$ is obvious. To prove the 
existence,  
one defines $d$ by the formula (4.4); 
the computations after (4.4) show that we get indeed a BV structure.   

Thus, we have proven 

{\bf 4.14. Theorem.} {\it The above constructions establish two 
inverse bijections between the set of mappings} (4.19) {\it 
satisfying} (4.20) {\it and the set of BV structures on $S(\CA)$.} \ $\btu$ 

This Theorem generalizes Theorems 4.10 and 4.3.

\bigskip
\centerline{\bf References}
\bigskip

[BB] A.~Beilinson, J.~Bernstein, A proof of Jantzen conjectures, in: 
I.M. Gelfand Seminar, S. Gelfand, S. Gindikin (Eds.), {\it Adv. in Sov. Math.} 
{\bf 16}, Part I, American Math. Society, Providence, Rhode Island, 1993, 
1-50. 

[BD1] A.~Beilinson, V.~Drinfeld, Chiral algebras I, Preprint.

[BD2] A.~Beilinson, V.~Drinfeld, Quantization of Hitchin's integrable system 
and Hecke eigensheaves, Preprint, 1997.  

[B] A.~Borel et al., Algebraic $D$-modules, Academic Press, Boston et al., 
1987. 

[D] V.~Drinfeld, Letter to the author, September 1988. 

[EGA III] A.~Grothendieck, J.~Dieudonn\'e, \'El\'ements de G\'eometrie Alg\'ebrique, III. 
\'Etude cohomologique des faisceaux coh\'erents (Premi\`ere Partie), 
{\it Publ. Math. IHES}, {\bf 11}, 1961.

[SGA 1] A.~Grothendieck, Rev\^etements \'etales et groupe fondamental, 
{\it Lect. Notes in Math.} {\bf 224}, Springer-Verlag, Berlin-Heidelberg-
New York, 1971. 

[H] V.~Hinich, Descent of Deligne groupoids, {\it IMRN}, 1997, No. 5, 
223-239. 

[KN] S.~Kobayashi, K.~Nomizu, Foundations of differential geometry, Volume I, 
Interscience Publishers, New York-London, 1963.   

[M] J.~Milnor, Singular points of complex hypersurfaces, Princeton University 
Press, Princeton, New Jersey, 1968. 

[R] Z.~Ran, Thickening Calabi-Yau spaces, in: Mirror Symmetry II, 
{\it AMS/IP Stud. Adv. Math.,} {\bf 1}, Amer. Math. Soc., Providence, 
Rhode Island, 1997, 393-400.  

[Sch] J.A.~Schouten, \"Uber Differentialkomitanten zweier kontravarianter 
Grossen, {\it Nederl. Acad. Wetensch. Proc.}, Ser. A, {\bf 43} (1940), 
449-452.  

[S] J.-P.~Serre, Alg\`ebre Locale Multiplicit\'es, {\it Lect. Notes in Math.}, 
{\bf 11}, Springer-Verlag, Berlin-Heidelberg-New York, 1975.

\bigskip

{\it Max-Planck-Institut f\"ur Mathematik,}

{\it Gottfried-Claren-Stra\ss e 26, 53225 Bonn}

e-mail: vadik\@mpim-bonn.mpg.de
 
\enddocument